\documentclass[a4paper,12pt]{article}
\usepackage{url}

\begin{document}

\begin{center}
\large\textbf{Allowable Changes of Variable for Recently Developed Improper Integral Definitions}

\large{Michael A. Blischke}
\end{center}

\begin{flushleft}
\textbf{Introduction}
\end{flushleft}

There are two basic types of improper integrals defined; integrals with an infinite limit, and integrals with finite limits where the integrand does not have a finite limiting value as the variable of integration approaches a particular ``critical'' value.  Improper integrals with more than one critical value, or with interior critical values, can be found as a sum of these two basic types.

Recently, in ~\cite{blischke:ALtDefInfLim} and ~\cite{blischke:ALtDefFinLim}, new definitions for the two types of improper integrals were introduced which allow for a greater range of functions to be integrated than the previous conventional definitions.  When using the new definitions, however, arbitrary changes of the variable of integration are not allowed, as they may lead to an integral with a different value.  Only simple changes of variable were previously shown to leave the value of the integral unchanged.  The purpose of this paper is to show sufficient conditions for performing changes of variable for each of the two types, without changing the value of the integral.

We use the definition from ~\cite{blischke:ALtDefInfLim} for an integral with an infinite upper limit and lower limit $a$,
\begin{equation}
\label{E001}
\textsf{Z} \hspace{-10.8 bp} \int^{\infty}_{a} f(x) \,dx \equiv \lim_{b\rightarrow\infty}\left\{\int^{b}_{a}f(x) \,dx + \int^{b+c}_{b} f(x) z(x-b) \,dx\right\}
\end{equation}
where $f(x)$ is the function to be integrated, and where $z(x)$ is a termination function.  The inclusion of the additional term (as compared with the conventional definition) inside the limit allows convergence to be rigorously shown for a greater range of functions, $f(x)$, than using the conventional definition.  The overstruck \sffamily Z \normalfont on integrals using the new definitions is included to distinguish them from integrals that exist using conventional definitions. 

For an improper integral with both limits finite, the general definitions are given in ~\cite{blischke:ALtDefFinLim}.  For the special case where the critical limit is the lower one, and is 0, with a noncritical upper limit $\beta$, the improper integral is defined as
\begin{equation}
\label{E010}
\textsf{Z} \hspace{-10.8 bp} \int_{0}^{\beta} g(u) \,du \equiv \lim_{\delta \to 0^+} \left\{\int_{0}^{\delta}  g(u) w(u/\delta) \,du + \int_{\delta}^{\beta} g(u) \,du \right\}.
\end{equation}

A full description, and many properties, of these integral definitions are found in ~\cite{blischke:ALtDefInfLim} and ~\cite{blischke:ALtDefFinLim}.

\begin{flushleft}
\textbf{Transformation between Infinite and Finite Limit Integrals}
\end{flushleft}

In ~\cite{blischke:ALtDefFinLim} a general change of variable is used to transform from one type to the other, given by
\begin{eqnarray}
\label{E025}
u \!\! &=& \!\! \psi(x) \\
x \!\! &=& \!\! \psi^{-1}(u) . \nonumber
\end{eqnarray}
It is required that $\psi(x)$ be finite for all $x > \psi^{-1}(\beta)$, that it be strictly monotonic with
\begin{equation}
\label{E026}
\psi'(x) < 0,
\end{equation}
and that
\begin{equation}
\label{E027}
\lim_{x \to \infty} \left\{ \psi(x) \right\} \equiv 0.
\end{equation}
These two properties imply that
\begin{equation}
\label{E027b}
\psi(x) > 0.
\end{equation}

In ~\cite{blischke:ALtDefFinLim} it is shown that when $\psi(x)$ satisfies (\ref{E026}) through (\ref{E027b}), the finite limit integral WRT $u$ and the corresponding infinite limit integral WRT $x$ are equal if they both exist.  We also have, for the particular case of
\begin{equation}
\label{E054}
\psi(x) = d \, e^{-\alpha x}
\end{equation}
with $\alpha > 0$, $d > 0$, that the two integrals are equivalent, with the existence of either one implying the existence of the other\footnote{This is shown in ~\cite{blischke:ALtDefFinLim} only for the particular case $d = 1$, but it is easily seen to hold for all $d > 0$.}.

\begin{flushleft}
\textbf{Change of Variable for Infinite Limit Integrals}
\end{flushleft}

Using the conditions for converting between the two types of integrals established in ~\cite{blischke:ALtDefFinLim}, we can easily determine sufficient conditions for changes of variable for each type, by converting to the other type and back.  We will first look at the infinite limit integral.

In the general case, we start with an integral WRT $x$, and use a change of variable to an integral WRT $y$,
\begin{equation}
\label{E200}
\textsf{Z} \hspace{-10.8 bp} \int_{a}^{\infty} f(x) \,dx = \textsf{Z} \hspace{-10.8 bp} \int_{a}^{\infty} f(Q(y)) \,dQ(y) = \textsf{Z} \hspace{-10.8 bp} \int_{a}^{\infty} f(Q(y)) Q'(y) \,dy
\end{equation}
where $x$ and $y$ are related as
\begin{eqnarray}
\label{E201}
y \!\! &=& \!\! P(x) \\
x \!\! &=& \!\! P^{-1}(y) \equiv Q(y) \nonumber
\end{eqnarray}
with
\begin{eqnarray}
\label{E202}
\lim_{x\to\infty} \left\{ P(x) \right\} = \infty \\
\lim_{y\to\infty} \left\{ Q(y) \right\} = \infty. \nonumber
\end{eqnarray}

From (\ref{E054}) we can always take
\begin{equation}
\label{E203}
u = e^{-y}.
\end{equation}

We then require $\psi(x)$ such that
\begin{equation}
\label{E204}
u = \psi(x)
\end{equation}
is a valid transform between finite and infinite limit integrals.  Thus
\begin{equation}
\label{E205}
e^{-y} = \psi(x) = e^{-P(x)}.
\end{equation}
Taking the derivative gives
\begin{equation}
\label{E206}
\psi'(x) = -P'(x) e^{-P(x)}.
\end{equation}

The criteria that must be satisfied are (\ref{E026}) through (\ref{E027b}).  The second is satisfied by (\ref{E202}), and the third is satisfied for all real $P(x)$.  So the transform
\begin{equation}
\label{E207}
y = P(x)
\end{equation}
is valid as long as (\ref{E026}) is satisfied, which just requires $P(x)$ to be monotonically increasing,
\begin{equation}
\label{E208}
P'(x) > 0 .
\end{equation}
Note that this ensures that $P^{-1}(y)$ in (\ref{E201}) exists.

In this case, with $P(x)$ real, and with (\ref{E202}) and (\ref{E208}) satisfied, performing the change of variable of integration (\ref{E200}) gives the same value as long as both integrals exist.  Some specializations of this are
\begin{equation}
\label{E209}
y = d \, x^r
\end{equation}
for $d > 0$, $r > 0$, and
\begin{equation}
\label{E210}
y = d \, e^{\alpha x}
\end{equation}
for $d > 0$, $\alpha > 0$.  The specialization (\ref{E209}) requires $a > 0$ in (\ref{E200}) if $r$ is not an odd integer.  

$\mbox{  }$

$\mbox{  }$

\begin{flushleft}
\textbf{Change of Variable for Finite Limit Integrals}
\end{flushleft}

Next, we'll look at changes of variable of integration for integrals with an improper limit at 0.  First, the specific change of variable
\begin{equation}
\label{E220}
u = d \, t^r  \mbox{   } \mbox{   } (r > 0, d > 0)
\end{equation}
is always a valid change of variable for the finite limit improper integral, giving
\begin{equation}
\label{E221}
\textsf{Z} \hspace{-10.8 bp} \int_{0}^{\beta} g(u) \,du = \textsf{Z} \hspace{-10.8 bp} \int_{a}^{\beta} g(t^r) \,d(t^r) = \textsf{Z} \hspace{-10.8 bp} \int_{0}^{\beta} g(t^r) r t^{(r-1)} \,dt.
\end{equation}

We can always convert from a finite limit to an infinite limit integral and back using the changes of variable
\begin{equation}
\label{E222}
u = d \, e^{-rx}
\end{equation}
and then
\begin{equation}
\label{E223}
t = e^{-x},
\end{equation}
so (\ref{E221}) is always valid.

In the general case, we have
\begin{equation}
\label{E230}
\textsf{Z} \hspace{-10.8 bp} \int_{0}^{\beta} g(u) \,du = \textsf{Z} \hspace{-10.8 bp} \int_{a}^{\beta} g(Q(t)) \,dQ(t) = \textsf{Z} \hspace{-10.8 bp} \int_{0}^{\beta} g(Q(t)) Q'(t) \,dt.
\end{equation}
where $u$ and $t$ are related as
\begin{eqnarray}
\label{E231}
t \!\! &=& \!\! P(u) \\
u \!\! &=& \!\! P^{-1}(t) \equiv Q(t) \nonumber
\end{eqnarray}
with
\begin{equation}
\label{E232}
P(u) > 0 \mbox{   } \mbox{   } (u>0),
\end{equation}
\begin{equation}
\label{E233}
P'(u) > 0 \mbox{   } \mbox{   } (u>0),
\end{equation}
and with
\begin{eqnarray}
\label{E234}
\lim_{u\to 0} \left\{ P(u) \right\} = 0 \\
\lim_{t\to 0} \left\{ Q(t) \right\} = 0. \nonumber
\end{eqnarray}
Again, (\ref{E233}) ensures that $P^{-1}(t)$ in (\ref{E231}) exists.

We can show that with the above requirements on $P(u)$, the change of variable from $u$ to $t$ is valid, so long as both integrals exist.

From (\ref{E054}) we can always take
\begin{equation}
\label{E235}
u = e^{-x}  \mbox{   or   } x = -\ln(u).
\end{equation}

We then find $\psi(x)$ such that
\begin{equation}
\label{E236}
t = \psi(x)
\end{equation}
is a valid transform between finite and infinite limit integrals, i.e. satisfies (\ref{E026}) through (\ref{E027b}).  We have
\begin{equation}
\label{E237}
\psi(x) = \psi\left( -\ln(u) \right) = P(u)
\end{equation}
so (\ref{E027b}) is satisfied.

Taking the derivative of $P(u)$ gives
\begin{equation}
\label{E238}
P'(u) = d/du \left( \psi(-ln(u)) \right)  = -\psi'\left( -\ln(u) \right) / u
\end{equation}
so (\ref{E026}) is satisfied.

Finally, using (\ref{E234}) and (\ref{E235}) we get that (\ref{E027}) is satisfied.  Thus the change of variable of integration (\ref{E230}) is valid as long as both integrals exist.

\begin{flushleft}
\textbf{Conclusion}
\end{flushleft}
Conditions have been presented for allowable changes of the variable of integration when using recently introduced definitions for finite limit and infinite limit improper integrals.  These changes of variable are valid when both the initial and final integrals exist.  In a restricted set of changes of variable, the existence of either integral guarantees existence of the other.  In general, however, the existence of one of the integrals does not guarantee the existence of the other.

\bibliographystyle{alpha}
\bibliography{ExtendedIntegration}

\end{document}